\documentclass [leqno,12pt]{article}
\usepackage[utf8]{inputenc} 
\usepackage[T1]{fontenc} 
\usepackage{amsmath}
\usepackage{amsfonts}
\usepackage{amsthm}
\usepackage{color}
\usepackage{graphicx}
\usepackage[section]{placeins}

\newlength{\oldparindent}
\newcommand{\un}{{\bf 1}}

\newcommand{\cL}{{\mathbb {L}}}

\newcommand{\bpf}{\begin{preuve}}
\newcommand{\epf}{ \end{preuve} \medskip}

\newcommand{\benum}{\begin{enumerate}}
\newcommand{\eenum}{\end{enumerate}}

\newcommand{\bitem}{\begin{itemize}}
\newcommand{\eitem}{\end{itemize}}

\newcommand{\brmq}{\begin{rmq}}
\newcommand{\ermq}{\end{rmq}}

\newcommand{\brmqs}{\begin{rmqs}}
\newcommand{\ermqs}{\end{rmqs}}

\newcommand{\bapp}{\begin{application}}
\newcommand{\eapp}{\end{application}}

\newcommand{\bapps}{\begin{applications}}
\newcommand{\eapps}{\end{applications}}

\newcommand{\bdefi}{\begin{definition}}
\newcommand{\edefi}{\end{definition}}

\newcommand{\beq}{\begin{equation}}
\newcommand{\eeq}{\end{equation}}

\def\bpm{\begin{pmatrix}}
\def\epm{\end{pmatrix}}

\newcommand{\bcas}{\begin{cases}}
\newcommand{\ecas}{\end{cases}}

\newcommand{\bex}{\begin{exemp}}
\newcommand{\eex}{\end{exemp}}

\newcommand{\bexs}{\begin{exemps}}
\newcommand{\eexs}{\end{exemps}}

\newcommand{\bprop}{\begin{proposition}}
\newcommand{\eprop}{\end{proposition}}

\newcommand{\bthm}{\begin{theoreme}}
\newcommand{\ethm}{\end{theoreme}}

\newcommand{\bcor}{\begin{corollaire}}
\newcommand{\ecor}{\end{corollaire}}

\newcommand{\blem}{\begin{lemme}}
\newcommand{\elem}{\end{lemme}}

\newcommand{\beqna}{\begin{eqnarray}}
\newcommand{\eeqna}{\end{eqnarray}}

\newcommand{\beqnas}{\begin{eqnarray*}}
\newcommand{\eeqnas}{\end{eqnarray*}}

\definecolor{green}{rgb}{0,.7,.2}
\definecolor{orange}{rgb}{0.9,.5,0}

\newcommand{\LL}{{\rm L}}

\def\tr{\textmd{trace}\,}

\def\det{{ \rm{det}}}  

\def\Id{{\rm{Id}}} 

\def\cC{{\mathcal C}}
\def\cD{{\mathcal D}}

\def\cG{{\mathcal  G}}
\def\cH{{\mathcal  H}}
\def\cI{{\mathcal  I}}
\def\cJ{{\mathcal  J}}

\def\cL{{\mathcal L }}
\def\cM{{\mathcal M }}
\def\cP{{\mathcal P }}
\def\cQ{{\mathcal Q }}

\def\cV{{\mathcal V}}

\def\bbB{{\mathbb{B}}}
\def\bbC{{\mathbb{C}}}

\def\bbG{{\mathbb{G}}}

\newcommand{\bbN}{{\mathbb {N}}}
\newcommand{\bbR}{{\mathbb {R}}}
\newcommand{\bbS}{{\mathbb {S}}} 

\newtheorem{theoreme}{Theorem}[section]
\newtheorem{lemme}[theoreme]{Lemma}
\newtheorem{definition}[theoreme]{Definition}
\newtheorem{proposition}[theoreme]{Proposition}
\newtheorem{corollaire}[theoreme]{Corollary}

\newenvironment{exemp}{\noindent{\bf Example. --- }}{\par}
\newenvironment{exemps}{\noindent{\bf Examples}\benum}{\eenum\par}
\newtheorem{rmq}[theoreme]{Remark}
\newtheorem{rmqs}[theoreme]{Remarks}

\newenvironment{preuve}{\noindent{\it Proof. --- }}
{\hfill\rule{1.3mm}{2mm}\par} 
\newenvironment{application}{\noindent{\bf Application. --- }}{\par}
\newenvironment{applications}{\noindent{\bf Applications. --- 
}\benum}{\eenum\par}

\theoremstyle{definition}

\author{D. Bakry \footnote{Institut de Math\'ematiques de Toulouse, Universit\'e Paul Sabatier, 118 
route de Narbonne, 31062 Toulouse, France}  }

\date{\today}

\date{\today}

\makeatletter \renewcommand{\@oddfoot}{\sl \small
 \hfil \thepage\hfil \today}
\renewcommand{\@oddhead}{\sl \small
 \hfil }
 \makeatother

\begin{document}

\title{Symmetric diffusions with polynomial eigenvectors
}
\maketitle

\abstract{We describe symmetric diffusion operators where the spectral decomposition is given through a family of orthogonal polynomials. In dimension one, this reduces to the case of Hermite, Laguerre and Jacobi polynomials. In higher dimension, some basic examples arise from compact Lie groups. We give a complete description of the bounded sets on which such operators may live. We then provide a classification of those sets when the polynomials are ordered according to their usual degrees. }

\par MSC : 60B15, 60B20, 60G99, 43A75,  14H50.

\par Key words : orthogonal polynomials, diffusion operators, random matrices.

\section{Introduction\label{sec.intro}}

Symmetric diffusion operators and their associated heat semigroups play a central rôle in the study of continuous Markov processes, and  also in differential geometry and partial differential equations. The analysis of the associated heat   or potential kernels have  been considered from many points of view, such as short and long time asymptotics,  upper and lower bounds, on  the diagonal  and away from it,  convergence to equilibrium, e.g.  All these topics had been deeply investigated along the past half century, see ~\cite{BGL, davies, varopouloscoulhonsaloffcoste92} for example. Unfortunately,  there are very few examples where computations are explicit. 

The spectral decomposition may provide one approach to the heat kernel study and the analysis of convergence to equilibrium, especially when the spectrum is discrete. Once again, there are very few models where this spectral decomposition is at hand, either  for the explicit expressions of the eigenvalues or the eigenvectors.  The aim of this survey is to present a family of models  where this spectral decomposition is completely explicit. Namely, we shall require the eigenvectors to be a polynomials in a finite number of variables. We are then dealing with orthogonal polynomials with respect to the reversible measure of this diffusion operator. Once again, orthogonal polynomial  have  been thoroughly investigated in many aspects, going back to the early works of Legendre, Chebyshev, Markov and Stieltjes, see~\cite{DX,Macdo1,szego75} e.g.

To be more precise, we shall consider some open connected  set $\Omega\subset \bbR^d$, with piecewise smooth boundary $\partial \Omega$ (say at least piecewise  $\cC^1$, and  may be empty), and some probability measure $\mu(dx)$ on $\Omega$, with smooth  positive density  measure $\rho(x)$ with respect to the Lebesgue measure  and such that polynomials are dense in $\cL^2(\mu)$. 

A diffusion operator $\LL$ on $\Omega$  (but we shall also consider such objects on smooth manifolds with no further comment)   is a  linear second order differential operator with no $0$-order terms, therefore written as
\beq\label{intro.def.diff}\LL(f) = \sum_{ij} g^{ij}(x) \partial^2_{ij} f + \sum_i b^i(x) \partial_i f,\eeq
 such that at every point $x\in \Omega$, the symmetric matrix $(g^{ij}(x))$ is non negative. For simplicity, we shall  assume here that this matrix is always non degenerate in $\Omega$ (it may, and will in general, be degenerate at the boundary $\partial \Omega$). This is an  ellipticity assumption, which will be in force throughout. We  will also assume that the coefficients $g^{ij}(x)$ and $b^i(x)$ are smooth. We are mainly interested in the case where $\LL$ is symmetric on $\cL^2(\mu)$ when restricted to smooth  functions compactly supported in $\Omega$. On the generator $\LL$, this translates into

\beq\label{intro.def.diff.sym}\LL (f) = \frac{1}{\rho} \sum_i \partial_i\Big(\rho \sum_j g^{ij} \partial_j f\Big), \eeq
as is readily seen using integration by parts in $\Omega$, see~\cite{BGL}.

We are also interested in the case when $\cL^2(\mu)$ admits a complete orthonormal basis $P_{q}(x), q\in \bbN$,  of polynomials such that $\LL (P_q) = -\lambda_q P_q$, for some real (indeed non negative) parameters $\lambda_q$.   This is equivalent to the fact that there exists   an increasing sequence $\cP_n$ of finite dimensional subspaces of the set $\cP$ of polynomials such that $\cup_n \cP_n= \cP$ and such that $\LL$ maps $\cP_n$ into itself.
 
When such happens, we have a spectral decomposition of $\LL$ in this basis, and, when a function $f\in \cL^2(\mu)$ is written as 
$f = \sum_q c_qP_q$, then $\LL(f) = \sum_q -\lambda_q c_q P_q$, such that any expression depending only on the spectral decomposition may be analyzed easily.

Our aim is to describe various families of such situations, which will be referred to as polynomial models. In dimension $d$,  many such models may be constructed with various techniques:  Lie groups, root systems, Hecke algebras, etc, see~\cite{ArCoVa01, Askey74,BracialiPerezPinar, Chihara78,DamMuYes13,DX,Heck,Heck1,Koorn1,Koorn2,Koorn3,Koorn4,KrallS, Macdo,Macdo1}, among many other possible references. Introducing weighted degrees, we shall analyse the situation where the operator maps  for any $k\in \bbN$ the space $\cP_k$ of polynomials with degree less $k$ into itself. For bounded sets $\Omega$ with regular boundaries, this leads to an algebraic description of the admissible boundaries for such sets. In dimension $2$ and for the usual degree, we give a complete description of all the admissible models (reducing to 11 families up to affine transformation).  We then present some other models, with other degrees, or in larger dimension, with no claim for exhaustivity.

This short survey is organized as follows. In Section~\ref{sec.dim1}, we present a brief tour of  the dimension 1 case, where the classical families of Hermite, Laguerre and Jacobi polynomials appear. We provide some geometric description for the Jacobi case (which will serve as a guide in higher dimension) together with various relations between those three families. In Section~\ref{sec.basics} we describe  some basic  notions  concerning the symmetric diffusion operators, introducing in particular the operator $\Gamma$ and the integration by parts formula. We also introduce the notion of image, or projection, of  operators, a key tool towards the construction of  polynomial models. In Section~\ref{sec.spheres}, we describe the Laplace operators on spheres, $SO(n)$ and $SU(n)$, and we provide  various  projections arising from these examples   leading to polynomials models, in particular for spectral measures in the Lie group cases. Section~\ref{sec.gal.case} describes the general case, with the introduction of the weighted degree models. In particular, when the set $\Omega$ is bounded, we provide the algebraic description of the sets on which such polynomial models may exist.  In particular, we show that the boundaries of those admissible sets $\Omega$ must lie in some algebraic variety, satisfying  algebraic  restrictions. The description of the sets lead then to the description of the measures and the associated operators. Section~\ref{sec.dim2} is a brief account of the complete classification for the ordinary degree of those bounded models in dimension 2. This  requires a precise analysis of the algebraic nature of the boundary which is only sketched in this paper. Section~\ref{sec.other.degrees} provides some examples of $2$ dimensional models with weighted degrees, and is far from exhaustive, since no complete description is valid at the moment. Section~\ref{sec.higher.dim} proposes some new ways (apart from tensorization) to construct higher dimensional models from low dimension ones.  Finally,  we give in Section~\ref{sec.pictures} the various pictures corresponding to the 2 dimensional models described in Section~\ref{sec.dim2}.

\section{Dimension 1 case : Jacobi, Laguerre and Hermite\label{sec.dim1}}
 In dimension one, given a probability measure $\mu$  for which the polynomials are dense in $\cL^2(\mu)$, there is, up to the choice of a sign, a unique family $(P_n)$ of polynomials with $\deg(P_n)=n$ and which is an orthonormal basis in $\cL^2(\mu)$. It is obtained through the Gram-Schmidt orthonormalization procedure   of the sequence $\{1,x, x^2, \cdots\}$. When does   such a sequence consists in  eigenvectors of some given diffusion operator of the form $L(f)(x)= a(x)\partial^2 f+ b(x)\partial f$ (a Sturm-Liouville operator)?  This had been described long ago   (see e.g.  \cite{BakM,Maze}), and reduces up to affine transformation to the classical cases of Jacobi, Hermite and Laguerre polynomials. 
 
 Those three families play a very important rôle in many branches of mathematics and ingeneering (mainly in statistics and probability for the  Hermite family and in fluid mechanics for the Jacobi one), and we refer to the huge literature on them for further information.  We briefly recall these  models.
 
\begin{enumerate}
\item $\cI= \bbR$. The measure $\mu$ is Gaussian centered:
$\mu(dx)=\frac{e^{-x^2/2}}{\sqrt{2\pi}}\,dx$. The associated polynomials are the Hermite polynomials $(H_n)$. They are eigenvectors of the 
 Ornstein--Uhlenbeck 
$$\cH=\frac{d^2}{dx^2}-x\frac{d}{dx}\,\quad  \cH(H_n)= -nH_n.$$
\item $\cI=\bbR_+^*$ The measure is the gamma measure   $\mu_a(dx)=C_ax^{a-1}e^{-x}\,dx$, $a>0$. The associated polynomials are the Laguerre polynomials $L_n^{(a)}$, and the associated operator is the  Laguerre operator 
$$\cL_a=x\frac{d^2}{dx^2}+ (a-x)\frac{d}{dx}, \quad \cL_a (L_n^{(a)})= -nL_n^{(a)}.$$
\item $\cI= (-1,1)$.  The measure is the beta measure $\mu_{a,b}(dx)=C_{a,b}(1-x)^{a-1}(1+x)^{b-1}\,dx$, $a,b>0$. The associated polynomials are the  Jacobi polynomials $(J_n^{(a,b)})$ and the associated Jacobi operator is 
$$\cJ_{a,b}=(1-x^2)\frac{d^2}{dx^2}- \big(a-b+(a+b)x\big)\frac{d}{dx}, \quad \cJ_{a,b} J_n^{(a,b)}= -n(n+a+b)J_n^{(a,b)}.$$

\end{enumerate}
 The Jacobi family contains the ultraspherical (or Gegenbauer) family (when $a=b$) with as particular cases the Legendre polynomials $a=b=0$, the Chebyshev first and second kind ($a=b= -1/2$ and $a=b=1/2$ respectively), which appear, after renormalization,  when writing $\cos(n\theta)= P_n(\cos\theta) $ (first kind) and $\sin(n\theta)= \sin(\theta)Q_n(\cos\theta)$ (second kind). The first two families (Hermite and Laguerre) appear as limits of the Jacobi case.  For example, when we chose $a=b=n/2$ and let then $n$ go to $\infty$, and scale the space variable $x$ into $x/\sqrt{n}$, the measure $\mu_{a,a}$ converges to the Gaussian measure, the Jacobi polynomials converge to the Hermite ones, and $\frac{2}{n} \cJ_{a,a}$ converges to $\cH$.   
  
  In the same way, the Laguerre case is obtained from the Jacobi one  fixing $b$, changing $x$ into $\frac{2x}{a}-1$, and letting $a$ go to infinity. Then $\mu_{a,b}$ converges to $\mu_b$, and $\frac{1}{a} \cJ_{a,b}$ converges to   $\cL_b$. 
  
  Also, when $a$ is a half-integer, the Laguerre operator may be seen as the image of the Ornstein--Uhlenbeck operator in dimension $d$. Indeed, as the product of one dimensional Ornstein--Uhlenbeck operators, the latter has generator 
  $H_d = \Delta -x.\nabla$. It's  reversible measure  is $e^{-|x|^2/2} dx/(2\pi)^{d/2}$,  it's  eigenvectors are the products $Q_{k_1}(x_1)\cdots Q_{k_d}(x_d)$, and the  associated process  $X_t= (X^1_t, \cdots, X^d_t)$, is formed of independent one dimensional Ornstein-Uhlenbeck processes. Then, if one considers $R(x)= |x|^2$, then one may observe that, for any smooth function $F: \bbR_+\mapsto \bbR$,
  $$\cH_d\big(F(R)\big) = 2\cL_a(F)(R),$$ where $a= d/2$.  In the  probabilist interpretation, this amounts to observe that if $X_t$ is a $d$-dimensional Ornstein--Uhlenbeck process, then $|X_{t/2}|^2$ is a Laguerre process with parameter $a= d/2$. 
  
  In the same way,  as we shall see below in Section~\ref{sec.spheres}, when $a=b= (d-1)/2$, $\cJ_{a,a}$ may be seen as the Laplace operator $\Delta_{S^{d-1}}$ on the unit sphere $\bbS^{d-1}$ in $\bbR^{d}$ acting on functions depending only on the first coordinate (or equivalently on functions invariant under the rotations leaving $(1, 0, \cdots, 0)$ invariant), and a  similar interpretation is valid for $J_{p/2, q/2}$ for  integers $p$ and $q$.  This interpretation comes from  Zernike and Brinkman \cite{ZerBrink} and Braaksma and Meulenbeld \cite{BraakMeul} (see also Koornwinder \cite{KoornWind}).   Jacobi polynomials  also play a central role in the analysis on compact Lie groups. Indeed, for $(a,b)$ taking the various values of  $(q/2, q/2)$, $((q-1)/2,1)$, $(q-1,2)$, $(2(q-1), 4)$ and $(4,8)$ the Jacobi operator $\cJ_{a,b}$ appears as the radial part of the Laplace-Beltrami (or Casimir) operator on the compact rank 1 symmetric spaces, that is spheres, real, complex and quaternionic projective spaces, and the special case of the projective Cayley plane (see Sherman \cite{Sherman}).

\section{Basics on symmetric diffusions\label{sec.basics}}

Diffusion operators are associated with diffusion processes (that is continuous Markov processes) through the requirement that, if $(X_t)$ is the Markov process with   associated generator $\LL$, then for any smooth function $f$, $f(X_t)-\int_0^t \LL(f)(X_s) ds$ is a local martingale, see~\cite{BGL}.    Here, we are mainly interested in such diffusion operators which are symmetric with respect to some probability measure $\mu$. In probabilistic terms, this amounts to require  that when the law of $X_0$ is $\mu$, then not only at any time the law of $X_t$ is still $\mu$, but also, for any time $t>0$, the law   of $(X_{t-s}, 0\leq s \leq t)$ is the same than the law of $(X_s, 0\leq s \leq t)$  (that is why this measure is often called the reversible measure).

 For any diffusion operator as such described in equation~\eqref{intro.def.diff}, we may define the carré du champ operator $\Gamma$ as the following bilinear application
 $$\Gamma(f,g)= \frac{1}{2} \Big( \LL(fg)-f\LL(g)-g\LL(f)\Big),$$ defined say for smooth functions defined on $\Omega\subset \bbR^d$.  From formula~\eqref{intro.def.diff.sym}, it is readily seen that 
 $$\Gamma(f,g)= \sum_{ij} g^{ij}\partial_if\partial_jg,$$
 such that $\Gamma(f,f)\geq 0$. As already mentioned, we restrict ourselves to the case where the matrix $(g^{ij})$ is everywhere positive definite. Then, the inverse matrix $(g_{ij})$ defines a Riemannian metric. By abuse of language, we shall refer to the matrix $(g^{ij})$ (or equivalently to the operator $\Gamma$) as the metric of the operator, although formally is should be called a co-metric. Then, Riemannian properties of this metric may carry important information on the models under study. Typically, most of the models described below in Section~\ref{sec.dim2} have constant curvature.
 
 The operator $\LL$ satisfies the following chain rule (or change of variable formula). For any $n$-uple $f= (f_1, \cdots, f_k)$ and any smooth function $\Phi : \bbR^k\mapsto \bbR$, 
 \beq\label{basics.chain.rule}
 \LL\big(\Phi(f)\big) = \sum_i \partial_i \Phi(f) \LL(f_i) + \sum_{ij} \partial^2_{ij} \Phi(f) \Gamma(f_i,f_j).
 \eeq
 This allows us to compute $\LL\big(\Phi(f)\big)$ as soon as we know $\LL(f_i)$ and $\Gamma(f_i,f_j)$. This is in particular the case for the coordinates $x^i : \Omega\mapsto \bbR$, where $\Gamma(x^i, x^j)= g^{ij}(x)$ and $\LL(x^i)= b^i(x)$, recovering then the form given in equation~\eqref{intro.def.diff}. But we may observe that the family $f= (f_i)$ in the previous formula does not need to be a coordinate system (that is a diffeomorphism  from $\Omega$ into some new set $\Omega_1$). There may be more function $(f_j)$ then required ($k> d$) or less ($k<d$). This remark will play an important rôle in the sequel.
 
 In general, when looking at such operators, one  considers first the  action of the operator  on smooth compactly supported function on $\Omega$. Since we want to work on polynomials, it is better to enlarge the set of functions we are working on, for example the set of $\cL^2(\mu)$ functions which are smooth and compactly supported  in a neighborhood of $\Omega$, referred to below just as "smooth functions".
 
 The operator is symmetric in $\cL^2(\mu)$ when, for any pair $(f,g)$ of smooth   functions 
 \beq\label{basics.symmetry} \int_\Omega \LL(f) \,g\, d\mu= \int_\Omega f\, \LL(g)\, d\mu.\eeq
 
 Usually, for this to be true, one should require one of the functions $f$ or $g$ to be compactly supported in $\Omega$, or ask for  some boundary conditions on $f$ and $g$, such as Dirichlet or Neuman.  However, in the case of polynomial models, the operator will be such that no boundary conditions will be required for equation~\eqref{basics.symmetry} to hold. More precisely, at any regular point of the boundary, we shall require the unit normal vector  to belong to the kernel of the matrix $(g^{ij})$. Under such assumption, the symmetry equation~\eqref{basics.symmetry} is satisfied whenever $f$ and $g$ are smooth functions.

 If we observe that $\LL(\un)=0$ (where $\un$ denotes the constant function), then applying~\eqref{basics.symmetry} to $\un$ shows that 
 $\int \LL(f)\, d\mu=0$, and therefore, applying the definition of $\Gamma$ and integrating over $\Omega$ (provided $fg\in \cL^2(\mu)$),  that
 \beq\label{basics.ipp}\int \LL(f)\, g \, d\mu = -\int \Gamma(f,g)\,d\mu= \int_\Omega f\, \LL g\, d\mu,\eeq which is called  the integration by parts formula.
 
 The fact that $\Gamma(f,f)\geq 0$ and the previous  formula~\eqref{basics.ipp} shows that, for any function~$f$,
 $\int_\Omega f \, \LL(f)\, d\mu \leq 0$, and therefore all eigenvalues are non positive.

The operator is entirely determined by the knowledge of $\Gamma$ and $\mu$, as is  obvious from formula \eqref{intro.def.diff.sym}, and the datum  $(\Omega, \Gamma, \mu)$ is  called a  Markov triple  in the language of \cite{BGL}, to which we refer for more details about the general description of such symmetric diffusion operators
 
 As already mentioned, we want to analyse those situations such that the set $\cP$ of polynomials is dense in $\cL^2(\mu)$ (and, being an algebra, all polynomials will be automatically in any $\cL^p(\mu)$ for any $p>1$), and such that  there exists some Hilbert basis $(P_n)$  of $\cL^2(\mu)$  with elements in $\cP$ such that  $\LL(P_n)= -\lambda_n P_n$. 
 Since we also require that any polynomial is a finite linear combination of the $P_n$', we see that the set $\cP_n= \{\sum_{0}^n \mu_k P_k\}$ is an increasing sequence of finite dimensional linear subspaces of $\cP$ such that $\LL : \cP_n\mapsto \cP_n$, and $\cup_n \cP_n= \cP$.
 
 Conversely, if there exists such an  increasing sequence $\cP_n$ of finite dimensional linear subspaces of $\cP$ such that  such that $\cup_n\cP_n= \cP$ satisfying $\LL : \cP_n\mapsto \cP_n$, then we may find a sequence $(P_n)$ which is an orthonormal basis of $\cL^2(\mu)$ and eigenvectors of $\LL$. Indeed, the restriction of $\LL$ to the finite dimensional subspace $\cP_n$ is symmetric when we provide  it with the Euclidean structure inherited from the $\cL^2(\mu)$ structure, and therefore may be diagonalized in  some orthonormal basis. Repeating this in any space $\cP_n$ provides the full sequence of polynomial orthogonal vectors.
 
 It may be worth to observe that when this happens, the set of polynomials is an algebra dense in $\cL^2(\mu)$ and stable under $\LL$ and the associated heat kernel $P_t= \exp(t\LL)$. When this happens, it is automatically dense in the $\cL^2(\mu)$ domain of $\LL$, for the domain topology,  and the set of polynomial will be therefore a core for our operator (see \cite{BGL} for more details). 
  
 From now on, we shall denote by $(P_n)$  such a sequence of eigenvectors, with $\LL(P_n)= -\lambda_n P_n$ (and we recall that $\lambda_n\geq 0$). Since $\LL(\un)=0$, we may always chose $P_0=\un$, and $\lambda_0=0$.  In general, this eigenvalue is simple, in particular in the elliptic case. Indeed, thanks to the integration by parts formula~\eqref{basics.ipp}, any function $f$ such that $\LL(f)=0$ satisfies 
 $\int \Gamma(f,f) d\mu=0$, from which $\Gamma(f,f)=0$. If ellipticity holds  (but also under much weaker requirements), then this implies that $f$ is constant.
 
 As mentioned in the introduction, one is often interested in the heat kernel $P_t$ associated with $\LL$, that is the linear operator $\exp(t\LL)$, defined through the fact that $\LL P_n = e^{-\lambda_n t} P_n$, or equivalently by the fact $F(x,t)= P_t(f)(x)$ satisfies the heat equation $\partial_t F= \LL_x(F)$, with $F(x,0)= f(x)$.
 
 This heat semigroup may be represented (at least at a formal level) as 
 $$P_t(f)(x) = \int_\Omega f(y) p_t(x,y) d\mu(y),$$ where the heat kernel may be written
 \beq\label{basics.rep.pt}p_t(x,y) = \sum_n e^{-\lambda_n t} P_n(x)P_n(y),\eeq
  provided that the series $\sum_n e^{-2\lambda_nt}$ is convergent (in practise and in all our models, this will always be satisfied). From this we see that a good knowledge on $\lambda_n$ and $P_n$ provides information of this heat kernel. However, it happens (thanks to the positivity of $\Gamma$) that $P_t$ preserves positivity, (and of course $P_t(\un)=1$), which is equivalent to the fact that $p_t(x,y)\mu(dy)$ is a probability measure for any $t>0$ and any $x\in \Omega$, in particular $p_t(x,y)\geq 0$. This is not  at all obvious from the representation~\eqref{basics.rep.pt}. Therefore, this representation~\eqref{basics.rep.pt} of $p_t(x,y)$ does not carry all the information about it.      
 
 It is worth to observe the following, which will be the basic tool for the construction of our polynomial models. Start from an operator  $\LL$ (defined on some manifold $\Omega_1$ or some open set in it), symmetric  under some probability measure $\mu$.  Assume that we may find  a set of functions $f= (f_i)_{i=1, \cdots k}$, that we consider as a function $f : \Omega\mapsto \Omega_1\subset \bbR^k$,  such that 
 $$L(f_i)= B^i(f), \Gamma(f_i,f_j)= G^{ij} (f).$$ 
 Then, for any smooth function $G : \Omega_1\mapsto \bbR$, and thanks to formula~\eqref{basics.chain.rule}, one has 
 $$\LL\big(G(f)\big)= \LL_1(G)(f),$$ where 
 \beq\label{image.diff}\LL_1 G= \sum_{ij} G^{ij} \partial^2_{ij} G+ \sum_i B^i\partial_i G.\eeq
 This new diffusion operator $\LL_1$ is said to be the image of $\LL$ under $f$.  In probabilistic terms, the image $Y_t= f(X_t)$ of the diffusion process $X_t$ associated  with $\LL$ under $f$ is still a diffusion process, and it's generator is $\LL_1$. Moreover, if 
 $\LL$ is symmetric with respect to $\mu$, $\LL_1$ is symmetric with respect of the image measure $\nu$ of $\mu$ under $f$. 
 
 This could (and will) be an efficient method  to determine the density $\rho_1$ of this image measure with respect to the Lebesgue measure, through the use of formula~\eqref{intro.def.diff.sym}.

\section{Examples : spheres, $SO(d)$, $SU(d)$\label{sec.spheres}}

We now describe a few natural examples  leading to       polynomial models. The first basic remark is that, given a symmetric diffusion operator $\LL$ described  as before by a Markov triple $(\Omega, \Gamma, \mu)$, it is enough to find a set $\{X_1, \cdots, X_k\}$  of functions such that for any $i= 1,\cdots, k$, $\LL(X_i)$ is a degree 1 polynomial in the variables $X_j$, and for any pair $(i,j), i,j= 1,\cdots, k$, $\Gamma(X_i,X_j)$ is a degree 2 polynomial in those variables. Indeed, if such happens, then the image $\LL_1$ of $\LL$ under $X= \{X_1, \cdots, X_k\}$ given in~\eqref{image.diff} is a symmetric diffusion operator  with reversible measure $\mu_1$, where $\mu_1$ is the image of $\mu$ under $X$.  Thanks to formula~\eqref{basics.chain.rule}, $\LL_1$  preserves for each $n\in \bbN$ the set  $\cP_n$ of polynomials  in the variables $(X_i)$ with total degree less than $n$.  One may then diagonalize $\LL_1$ in $\cP_n$, and this leads to the construction of a $\cL^2(\mu_1)$ orthonormal basis formed with polynomials.

Moreover,  given polynomial models, we may consider product models, which are again polynomial models, and from them consider projections. Indeed,  given two  polynomial models   described by  triples $(\Omega_i, \Gamma_i, \mu_i)$ as in  Section~\ref{sec.basics}, we may introduce  on $\Omega_1\times \Omega_2$ the  measure $\mu_1\otimes \mu_2$,  and the sum operator $\Gamma_1\oplus \Gamma_2$, acting on functions $f(x,y)$ and $h(x,y)$ as $$\big(\Gamma_1\otimes \Id\oplus\Id\otimes \Gamma_2\big)(f,h)(x,y)= \sum_{ij}g_1^{ij} \partial_{x_i}f \partial_{x_j}h+  \sum_{kl}g_2^{kl}(y) \partial_{y_k}f \partial_{y_l}h.$$

Once we have some polynomial models $(\Omega_1, \Gamma_1, \mu_1)$ and  $(\Omega_2, \Gamma_2, \mu_2)$, then 
$(\Omega_1\times \Omega_2, \Gamma_1\oplus\Gamma_2, \mu_1\otimes\mu_2)$ is again a polynomial model. At the level of processes, if $(X^i_t)$ is associated with $(\Omega_i, \Gamma_i, \mu_i)$, and are chosen independent, then the process $(X^1_t, X^2_t)$ is associated with this product. This kind of tensorisation procedure constructs easily new dimensional models in higher dimension once we have some in low dimension.  The  polynomials $R_{i,j}(x,y), x\in \Omega_1, ~ y\in \Omega_2$ associated with the product  are then just the tensor  products $P_i(x)Q_j(y)$  of the polynomials associated with each of the components.

Moreover, one may consider quotients in these products to construct more polynomial models, as we did to pass from the one dimensional Ornstein-Uhlenbeck operator to the Laguerre operator.

The easiest example to start with  is the Laplace operator $\Delta_\bbS$ on  the unit sphere  $\bbS^{d-1}\subset \bbR^d$. This operator may be naively described as follows : considering some  smooth function $f$ on $\bbS^{d-1}$, we  extend it in a neighborhood of $\bbS^{d-1}$ into a function which is independent of the norm, that is $\hat f(x)= f(\frac{x}{\|x\|})$, where $\|x\|^2= \sum_i x_i^2$, for $x=(x_1, \cdots, x_d)\in \bbR^d$. Then, we consider $\Delta(\hat f)$, where $\Delta $ is the usual Laplace operator in $\bbR^d$, and restrict $\Delta(\hat f)$ on $\bbS^d$. This is $\Delta_\bbS (f)$.

An easy exercise shows that, for the functions $x_i$ which are the restriction to $\bbS^{d-1}$ of the usual coordinates in $\bbR^d$, then
$$\Delta_\bbS(x_i)= -(d-1)x_i, ~\Gamma_\bbS(x_i,x_j) = \delta_{ij}-x_ix_j.$$
The uniform measure on the sphere (that is the unique probability  measure which is invariant under rotations) is the reversible measure for $\Delta_\bbS$. 
The system of functions $(x_1, \cdots, x_d)$ is not a coordinate system, since those functions are linked by the relation $\sum_i x_i^2= 1$. For example, one sees from the previous formulae that $\Delta_\bbS (\|x\|^2)=\Gamma(\|x\|^2, \|x\|^2)=0$ on $\bbS^{d-1}$ (an good  way to check the right value for $\lambda$ when imposing $\LL  (x_i)= -\lambda x_i$ with the same $\Gamma$ operator).

 But the system $(x_1, \cdots, x_{d-1}) $ is a system of coordinates for say the upper half sphere. We may observe that the operator indeed projects onto $(x_1, \cdots, x_{d-1})$ into an elliptic operator in the unit ball $\bbB^{d-1}=\{\|x\|< 1\}$, with exactly the same relations for $\LL(x_i)$ and $\Gamma(x_i,x_j)$. The image operator (that is the Laplace operator in this system of coordinates)  is then
$$\LL= \sum_{ij} (\delta_{ij} -x_ix_j) \partial^2_{ij} -(d-1)\sum_i x_i \partial_i,$$ and from the formula~\eqref{intro.def.diff.sym}, is is easy to determine the density measure  (up to a multiplicative constant), which is
$\rho(x)= (1-\|x\|^2)^{-1/2}$, which happens here to be $\det(g^{ij})^{-1/2}$ (as is always the case for Laplace operators).  We see then that this provides an easy way to compute the density of the uniform measure on the sphere under this map $\bbS^{d-1}\mapsto \bbB^{d-1}$, which projects the upper half sphere onto the unit ball.

One may also observe that if $M= (M_{ij})$ is a matrix (with fixed coefficients), and if $y_i= \sum_j M_{ij} x_j$, then
$\Delta(y_i)= -(d-1)y_i$ and $\Gamma(y_i,y_j)= (MM^t)_{ij}-y_iy_j$. Then, when $M$ is orthogonal, the image measure of $\Delta_\bbS$ under $x\mapsto Mx$ is $\Delta_\bbS$ itself, which tells us that the Laplace operator is invariant under  orthogonal transformations.

We may also consider the projection $\pi$  from $\bbS^{d-1}$ to $\bbB_p$ for $p<d-1$ : $\pi :(x_1, \cdots , x_d)\mapsto   (x_1, \cdots, x_p)$, which provides the same operator as before, except that now we are working on the unit ball $\bbB^p\subset \bbR^p$ and $L(x_i)= -(d-1)x_i$, where the parameter  $d$ is no longer correlated with the dimension $p$ of the ball. We may as well consider the generic operator  $\LL_\lambda$ on the unit ball in $\bbR^p$ with 
$\Gamma(x_i,x_j)= \delta_{ij} -x_ix_j$ and $\LL(x_i)= -\lambda x_i$, where $\lambda>p-1$. Is is readily checked that it has symmetric measure density $C_{p,d} (1-\|x\|^2)^{\frac{\lambda-p-1}{2}}$. As a consequence, the image measure of the sphere $\bbS^{d-1}$ onto the unit ball through this projection has density  $C_{d,p}(1-\|x\|^2)^{(d-p-2)/2}$.  It is worth to observe that when $\lambda$ converges to $p-1$, the measure converges to the uniform measure on the boundary of $\bbB^p$, that is $\bbS^{p-1}$, and the operator converges to the operator on $\bbS^{p-1}$. 

When we chose $p=1$, we  recover the symmetric Jacobi polynomial model in dimension $1$ with parameters $a=b= (d-1)/2$.

For these operators, we see that, in terms of the variables $(x_i)$, $g^{ij}$ are polynomials with total degree $2$ and $L(x_i)$ are  degree 1. Therefore, in view of the chain rule formula~\eqref{basics.chain.rule},  we see that the operator $\LL$ such defined maps the space $\cP_n$ of polynomials with total degree less than $n$ into itself, and this provides  a first family of polynomial models.

One may also still consider the unit sphere in $\bbR^d$, and chose integers such that  $p_1+ \cdots+ p_k=d$. Then, setting $P_0=0, P_i= p_1+ \cdots +p_i$, consider the functions 
$X_i= \sum_{P_{i-1}+1}^{P_{i}} x_j^2$, for $i= 1, \cdots k-1$. The image of the sphere under this application is the simplex $\cD_{k-1}= \{X_i\geq 0, \sum_1^{k-1} X_i \leq 1\}$. 
We have $\Delta_\bbS(X_i)= 2(p_i-dX_i)$, and $\Gamma(X_i,X_j)= 4X_i(\delta_{ij}-X_iX_j)$.
The operator $\Delta_\bbS$ such projects on the simplex,  with $$G^{ij} = 4X_i(\delta_{ij}-X_j), ~B^i= 2(p_i-dX_i),$$ and provides again a polynomial model on it. The reversible measure is easily seen to have density
$C\prod_1^{k-1} X_i^{r_i}(1-\sum_1^{k-1}X_i)^{r_k}$, with $r_i= \frac{p_i-2}{2}, ~ri= 1, \cdots, k$, which a Dirichlet measure on the simplex $\cD_{k-1}$. This measure is then seen as the image of the uniform measure on the sphere through this application $X : \bbS^{d-1}\mapsto \cD_{k-1}$. The general Dirichlet density measure $\rho_{a_1, \cdots, a_k}(X)= X_1^{a_1}\cdots X_{k-1}^{a_{k-1}} (1-\sum_i X_i)^{a_k}$ also produces (with the same $\Gamma$ operator)  a new family of polynomial models.  We may play a around  with the Dirichlet measures for half integer parameters.  For example, look at the image   in a Dirichlet $\rho_{a_1, \cdots, a_k}$ measure  for the image through $(X_1,X_2, \cdots, X_{k-1})\mapsto (X_1+X_2, X_3, \cdots, X_{k-1})$ produces a Dirichlet measure with parameters $\rho_{a_1+a_2, a_3, \cdots, a_k}$, which is obvious from the sphere interpretation, and extends easily to the general class. The same procedure is still true at the level of the operators (or the associated stochastic processes).

Once again, when considering the case $k=2$, we get an operator on $[0,1]$. Changing $X$ into $2X-1$, we get an operator on $[-1,1]$, which, up to the scaling factor 4,  is the dissymetric Jacobi model with parameters $a= r_1/2, b= (d-r_1)/2$. 

There are many other ways to produce polynomial models from spheres, (and we shall provide some later in dimension $2$, see Section~\ref{sec.other.degrees}. But we want to show another very general way to construct  some. Let us consider some semi-simple compact Lie group $\bbG$ (such as $SO(n), SU(n), Sp(n)$, etc). On those groups, there exist a unique invariant probability measure (that is invariant under $x\mapsto gx$ and $x\mapsto xg$, for any $g\in \bbG$). This is the Haar measure. There  also exists a unique (up to a scaling constant) diffusion operator  $\Delta_\bbG$ which is also invariant under left and right action : this means that if, for a smooth  function $f : \bbG\mapsto \bbR$, one defines for any $g\in \bbG$,  $L_g(f)(x)= f(gx)$, then $\Delta_\bbG (L_g(f))= L_g(\Delta_\bbG f)$, and the same is true for the right action $R_g(f)(x)= f(xg)$.  This operator is called the Laplace (or Casimir ) operator on the group $\bbG$.  Assume then that the Lie group $\bbG$ is represented as a group of matrices, as  is the case for the natural presentation of the natural above mentioned  Lie groups. The Lie algebra $\cG$ of $\bbG$ is the tangent space at the origin, and to any element $A$ of the Lie algebra, that is a matrix $(A_{ij})$, we may associate some vector field $R_A$ on the group through $X_A(f)(g) = \partial_t f(ge^{tA}) _{t=0} $.  If we write  $g= (g_{ij})$ and consider a function $f(g_{ij}) : \bbG\mapsto \bbR$, then 
$$X_A(g)= \sum_{ijk}  g_{ik}A_{kj} \partial_{g_{ij}} f,$$ and therefore $X_A$ preserves the set $\cP_n$ of polynomials with total degree $n$ in the variables $(g_{ij})$.
Now, the  Casimir operator may be written as $\sum_i X_{A_i}^2$, where the $A_i$ form an orthonormal basis for some natural quadratic form on the Lie algebra $\cG$ called the Killing form. This operator also preserves the set $\cP_n$. Unfortunately, those "coordinates" $g_{ij}$ are in general linked by algebraic relations, and may not serve as a true coordinate system on the group. However, we may then describe the operator $\Delta_\bbG$ through it's action on those functions $g_{ij} : \bbG\mapsto \bbR$.

Without further detail, consider  the group $SO(n)$ of orthogonal matrices with determinant 1. Let  $m_{ij}$ be the entries of a matrix in $SO(n)$, considered as functions on the group. We have
$$\Delta_{SO(n)}(m_{ij})=-(n-1) m_{ij}, ~\Gamma_{SO(n)}( m_{kl},m_{qp})= \delta_{kq}\delta_{lp}- m_{kp}m_{ql}.$$

We now show some projections of this operator. Let $\cM(p,q)$ be  the space of $p\times q$ matrices. Select   $p$ lines and $q$ columns in the matrix $g\in SO(n)$, (say the first ones), and consider the map $\pi : SO(n)\mapsto \cM(p,q)$ which to $M\in SO(n)$ associates  the extracted matrix $N= (m_{ij}), 1\leq i \leq p, 1\leq j \leq q$. From the form of $\Delta_{SO(n)} (m_{ij})$ and $\Gamma_{SO(n)}(m_{ij}, m_{kl})$, it is clear  that   the operator projects on $\cM(p,q)$  through $\pi$. It may happen (whenever $p+q> n$)  that  the image is contained in a sub-manifold (indeed an algebraic variety) of $\cM(p,q)$.  But we have nevertheless a new diffusion operator on this image,  and the associated process is known as the matrix  Jacobi process. It is worth to observe that if $p=n$ and $q=1$, this is nothing else than the spherical operator $\Delta_\bbS$ in $\bbS^{n-1}$. In general, whenever $p+q\leq n$, this process lives of the symmetric domain $\{NN^*\leq \Id\}$, and  has a reversible measure with density $\rho$ with respect to the Lebesgue measure which is $\det(\Id-NN^*)^{(n-p-q-1)/2}$, which is easily seen from formula~\eqref{intro.def.diff.sym}. We may also now fix $p$ and $q$ and consider $n$ as a parameter,  and we obtain a family of polynomial processes  on this symmetric domain a long  as $p+q<n+1$.
 
One may play another game and consider the image of the operator on the spectrum. More precisely, given the associated process $X_t\in SO(n)$, one looks at the process  obtained on the eigenvalues of $X_t$ (that is the spectral measure of $X_t$). This process is again a diffusion process, for which we shall  compute the generator.
To  analyze the spectral measure (that is the eigenvalues up to permutations of them), the best is to look at the characteristic polynomial $P(X)= \det(M-X\Id)= \sum_{i=1}^n a_i X^i$. Then we wan to compute $\Delta_{SO(n)}(a_i)$ and $\Gamma_{SO(n)}(a_i,a_j)$.

For a generic matrix $M= (m_{ij})$, Cramer's formulae tells us that, on the set where $M$ is invertible,  $\partial_{m_{ij}} \log (\det(M))=M^{-1}_{ji}$ and 
$\partial^2_{m_{ij}, m_{kl}} \log(\det (M))= -M^{-1}_{jk}M^{-1}_{li}$.

From this,  and using the chain rule formula, we get that 
\beqnas \Delta_{SO(n)} \log P(X)=&& -(n-1)\tr(M(M-X\Id)^{-1}\\&&- \tr\Big((M-X\Id)^{-1}(M^t-X\Id)^{-1}\big)+ \Big(\tr M(M-X\Id)^{-1}\Big)^2.
\eeqnas
and 
\beqnas \Gamma\big(\log (P(X)), \log (P(Y))\big)= &&\tr\big((M-X\Id)^{-1}(M^t-Y\Id)^{-1}\big)\\&&- \tr\big(M^2(M-X\Id)^{-1}(M-Y\Id)^{-1}\big).\eeqnas

But 
$$\bcas \tr(M(M-X\Id)^{-1}= n-X\frac{P'}{P}(X), \\
\tr\big((M-X\Id)^{-1}(M^t-Y\Id)^{-1}\big)= \frac{1}{1-XY}\Big(\frac{1}{Y}\frac{P'}{P}(\frac{1}{Y})-X\frac{P'}{P}(X)\Big)\\
 \tr\big(M^2(M-X\Id)^{-1}(M-Y\Id)^{-1}\big)= n+ \frac{1}{X-Y}\Big(X^2\frac{P'}{P}(X)-Y^2\frac{P'}{P}(Y)\Big).
\ecas$$
One may use the fact that $M\in SO(n)$ to see that $P(\frac{1}{X})= (-X)^{-n}P(X)$, so that 
$\frac{1}{Y}\frac{P'}{P}(\frac{1}{Y})=n-Y\frac{P'}{P}(Y)$.

In the end, we see that 
$$\Delta_{SO(n)} (P)= -(d-1) XP' + X^2P'' ,$$
$$ \Gamma_{SO(n)}\big(P(X),P(Y)\big)= \frac{XY}{1- XY}\Big(n P(X)P(Y)+ \frac{(1-X^2)P(Y)P'(X)-(1-Y^2)P(X)P'(Y)}{X-Y}\Big).$$
Since $\Delta_{SO(n)} P(X) = \sum_i \Delta_{SO(n)} (a_i) X^i$ and 
$$\Gamma_{SO(n)}\big(P(X),P(Y)\big)=\sum_{ij} \Gamma_{SO(n)}(a_i,a_j) X^iY^j,$$ we see form the action  of $\Delta_{SO(n)}$ and $\Gamma_{SO(n)}$ on $P(X)$  that, in terms of the variables $(a_k)$, $\Delta_{SO(n)}(a_i)$ are   degree one   polynomials  and  $\Gamma_{SO(n)}(a_i,a_j)$ are  degree two. Therefore, from the same argument as before, the operator $\Delta_{SO(n)}$ projects on it's spectrum into a polynomial model in the variables $a_i$.

The same is true  (with similar computations) for the spectra of $NN^*$, where $N$ is the extracted matrix in $\cM(p,q)$ described above (corresponding to the matrix Jacobi process), and for many more models. 

Similarly, one may also look  at the special unitary group $SU(n)$, where the coordinates $(z_{ij}= x_{ij}+ {\bf i} y_{ij})$  are the entries of the matrix.  Using complex coordinates, one has then to consider $\Delta_{SU(n)}(z_{ij})$ and $\Gamma(z_{ij}, z_{kl})$ and $\Gamma(z_{ij}, \bar z_{kl})$ in order to recover the various quantities corresponding to the variables $x_{ij}$ and $y_{ij}$ (using the linearity of $\Delta$ and the bilinearity of $\Gamma$). We have (up to some normalization)
$$\bcas \Delta_{SU(n)}(z_{ij})= -(n^2-1)z_{ij},\\
 \Gamma_{SU(n)}(z_{ij},z_{kl})= z_{ij}z_{kl}- nz_{il}z_{kj}, \\
\Gamma_{SU(n)}(z_{ij},\bar z_{kl})= n\delta_{ik}\delta_{jl} -z_{ij}\bar z_{kl})
\ecas
$$
The same remark as before applies about the extracted matrices, and also, with the same method,  we get for the characteristic polynomial
\beq\label{spect.su3}\bcas \Delta_{SU(n)}(P)= -(n^2-1) X P' + (n+1)X^2 P'',\\
 \Gamma_{SU(n)}(P(X), P(Y))= XY\Big(P'(X)P'(Y)+ \frac{n}{X-Y}\big(P'(X)P(Y)-P'(Y)P(X)\big)\Big),\\
 \Gamma_{SU(n)}(P(X), \bar P(Y))=\frac{1}{1-X Y}\Big(nP(X)\bar P(Y)- \bar Y \bar P'( Y)P(X)- XP'(X)\bar P(Y)\Big).
 \ecas
 \eeq
 
 Once again, the Casimir operator on $SU(n)$ projects onto a polynomial model in the variables of the characteristic polynomial.

\section{The general case\label{sec.gal.case}}

As we briefly showed in Section~\ref{sec.spheres}, there are  many models for orthogonal polynomials and they are quite hard to describe in an exhaustive way.  We  propose in this Section  a more systematic approach. Recall that we are considering probability measures $\mu$  on $\bbR^d$, on some open connected set   $\Omega\subset \bbR^d$ for which the set $\cP$ of polynomials are dense in $\cL^2(\mu)$. Recall that it is enough for this to hold that there exists some $\epsilon>0$ such that $\int e^{\epsilon\|x\|} d\mu <\infty$.

The first thing to do is to describe  some ordering of the polynomials. For this, we chose  a sequence $a=(a_1, \cdots, a_d)$ of positive integers  and say that the degree of a monomial $x_1^{p_1}x_2^{p_2} \cdots x_d^{p_d}$ is $a_1p_1+ \cdots a_dp_d$. Then, the degree $\deg_a(P)$ of a polynomial $P(x_1, \cdots, x_d)$ is the maximum value of the degree of it's monomials (such a degree is usually referred as a skew  degree). Then, the space $\cP_n$ of polynomials with 
$\deg_a(P)\leq n$ is finite dimensional. Moreover, $\cP_n\subset \cP_{n+1}$, and $\cup_n \cP_n= \cP$. 

To chose a  sequence of orthogonal polynomials $(P_k)$ for $\mu$, we chose at each step $n$ some orthonormal basis in the orthogonal complement $\cH_n$ of $\cP_{n-1}$ in $\cP_n$.   This space in general has a large dimension (increasing with $n$), and there is therefore not an unique choice of such an orthonormal basis. 

We are then looking for diffusion differential operators $\LL$  (with associated $\Gamma$ operator) on $\Omega$, such that $\LL$ admits such a sequence $(P_n)$ as eigenvectors, with real eigenvalues . The operator $\LL$ will be then automatically  essentially self-adjoint on $\cP$. 
The first observation is that for each $n$, $\LL$ maps $\cP_n$ into itself.

In particular, for each coordinate, $\LL(x_i)\in \cP_{a_i}$ and for each pair of coordinates $(x_i,x_j)$, $\Gamma(x_i,x_j)\in \cP_{a_i+a_j}$. 
Then, writing $\LL= \sum_{ij} g^{ij} \partial_{ij}+ \sum_i b^i\partial_i$,  we see that 
\beq\label{cond.deg}g^{ij}\in\cP_{a_i+a_j}, ~b^i\in \cP_{a_i}\eeq
Moreover, under the conditions~\eqref{cond.deg}, we see from the chain rule formula~\eqref{basics.chain.rule} that for each $n\in \bbN$, $\LL$ maps $\cP_n$ into itself. Provided it is symmetric  on $\cP_n$ for the Euclidean structure induced from $\cL^2(\mu)$, we will then be able  to derive an orthonormal basis formed with eigenvectors for it.

Once we have a polynomial model with a given choice of degrees $(a_1, \cdots, a_d)$, (say, in the variables $(x_i, i= 1, \cdots, d)$,  and as soon as one may find polynomials $X_1, \cdots, X_k$ in the variables $x_i$ with $\deg_a(X_i)= b_i$,  as soon as $\LL(X_i)$ and $\Gamma(X_i,X_j)$ are polynomials in those variables $X_i$, then we get a new model in the variables $(X_i)$ (provided however that the ellipticity requirement is satisfied), with new degrees $b=(b_1, \cdots, b_k)$. Indeed, from the chain rule~\eqref{basics.chain.rule}, one sees that the image  operator $\LL_1$ maps the set $\cQ_n$ of polynomials in the variables $(X_i)$  with degree $\deg_b\leq n$ into itself. 

The next task is then  to describe the sets $\Omega$ on which such a choice for  $\LL$ and $\mu$ is possible. For that, we restrict our attention for those $\Omega\subset \bbR^d$ which are bounded with piecewise $\cC^1$ boundaries. We call those sets admissible sets.

Then, we have the main important characterization
\bthm\label{main.thm.gal}~

\benum
\item
If $\Omega$ is an admissible set, then $\partial\Omega$ is included into an algebraic variety, with $\deg_a \leq 2\sum_i a_i$.
\item Let $Q$ be the reduced equation of $\partial_\Omega$. $\Omega$ is admissible if and only if there exist some $g^{ij}\in \cP_{a_i+a_j}$ and $L_i\in \cP_{a_i}$  such that 
\beq\label{cond.ipp.bord} \forall i= 1, \cdots, d, \sum_{j} g^{ij} \partial_j Q= L_i Q, ~(g^{ij}) \hbox{ with  non negative in $\Omega$}.\eeq
When such happens, $Q$ divides $\det(g)$.

\item 
Let $Q= Q_1\cdots Q_k$ the decomposition of $Q$ into irreducible factors.
If~\eqref{cond.ipp.bord} is satisfied, then any  function $\rho(x)= Q_1^{r_1}\cdots Q_k^{r_k}$ is an admissible density measure for the measure $\mu$, provided $\int_\Omega \rho(x) dx < \infty$. When $Q= \det(g)$, then there are no other measures.
\item For any solution $(g^{ij})$ of~\eqref{cond.ipp.bord}, and any $\mu$ as before, setting $\Gamma(f,g)= \sum_{ij} g^{ij}\partial_i\partial_jg$, then the triple $(\Omega, \Gamma, \mu)$ is an admissible solution for the operator $\LL$.
\eenum

\ethm
 \brmq Observe that equation~\eqref{cond.ipp.bord} may be rewritten as $\Gamma(\log Q, x_i)= L_i$.
 
 \ermq
\bpf We shall not give the full details of the proof here, and just describe   the main ideas. 

 Suppose we have a polynomial model with coefficients $g^{ij}$, $b^i$ on $\Omega$, with polynomial functions $g^{ij}$ and $b^i$ satisfying the above requirements on their degrees.

The first thing to observe is that if $\LL$ is diagonalizable of $\cP_n$ for each $n\in \bbN$, then for each polynomial pair $(P,Q)$
\beq\label{ipp.polyn}\int_\Omega \LL(P)\, Q\, d\mu= \int_\Omega P\, \LL (Q)\, d\mu.\eeq 
 This extends easily to any pair $(f,g)$  of smooth functions compactly supported in $\Omega$, so that the description~\eqref{intro.def.diff.sym} holds true. Moreover, $\Omega$ being bounded,  and the coefficients $g^{ij}$ and $b^i$ being polynomials, formula~\eqref{ipp.polyn} extends further to every pair $(f,g)$ of smooth functions, not necessarily with support in $\Omega$. Using Stockes formula, (and the regularity of the boundary of $\Omega$), this imposes that, for any pair of smooth function $(f,h)$, $\int_{\partial \Omega}\sum_{ij} f\partial_ih  g^{ij}n_j dx =0$, where $(n_j)$ is the normal tangent  vector at the boundary $\partial\Omega$. Therefore, this implies that, for any $i$, $\sum_jg^{ij}n_j=0$ on the boundary, so that $(n_j)$ is in the kernel of $(g)$ at the boundary. This implies in turn that the boundary lies inside the algebraic set $\{\det(g)=0\}$. 

Therefore, $\partial\Omega$ is included in some algebraic variety. For any regular point $x\in \partial\Omega$ consider an irreducible polynomial $Q_1$ such that, in a neighborhood $\cV$ of $x$, the boundary is included in   $\{Q_1=0\}$. Then, $(n_j)$ is parallel to $\partial_j Q_1$, so that 
$\sum_j g^{ij} \partial_j Q_1=0$ on $\cV\cap \{Q_1=0\}$. From Hilbert's Nullstellensatz, $\sum_j g^{ij} \partial_j Q_1= L_iQ_1$, for some polynomial $L_i$.

 This being valid for any polynomial $Q_1$ appearing in the reduced equation of $\partial\Omega$, this is still true for the reduced equation itself (and in fact the two assertions are equivalent). 

For a given polynomial $Q$, if equation~\eqref{cond.ipp.bord}  admits a non trivial solution $(g^{ij})$, then $\partial_iQ$ is in the kernel of $(g^{ij})$ at every regular point of $\{Q=0\}$. Then, $\det(g)$ vanishes at that point. $Q$ being reduced, then $Q$ is a factor of $\det(g)$.

Now, the link between $b^i= \LL(x_i)$ and $\sum_ g^{ij}\partial_j \log \rho$ given in~\eqref{intro.def.diff.sym} shows that, in order for $b^i$ to be a polynomial with degree less than $a_i$, it is enough (and in fact equivalent) to have $\sum_j g^{ij}\partial_j \log \rho=A_i$, for some polynomial $A_i$ with degree less than $a_i$. 
But comparing with equation~\eqref{cond.ipp.bord}  shows that it is satisfied for $Q_i^{r_i}$ for any factor $Q_i$ of $Q$ and any parameter $r_i$.
Then, all the condition are satisfied and the model $(\Omega, \Gamma, \mu)$ is a polynomial model.

\epf

One sees that indeed the problem of determining polynomial models relies entirely on the study of the boundary $\partial \Omega$, at least as far as bounded sets $\Omega$ are considered. Given any algebraic curve, and a fixed choice of degrees $(a_1, \cdots, a_k)$ it is an easy task to decide if this curve is a candidate to be (or to belong to) the boundary of some set $\Omega$ on which there exist a polynomial model : equation~\eqref{cond.ipp.bord} must have a non trivial solution. This equation is a linear system of equations in the coefficients of the polynomials $g^{ij}$ and $L_i$, however  in general  with much more equations than variables.  

Moreover, as soon as one has a model on $\Omega$, there there exist as we already saw  many other models on the same set, with the same $(g^{ij})$,  with measures described with a finite number of parameters, depending on the number of irreducible components in  the reduced equation of $\partial \Omega$.
 
The solutions of equation~\eqref{cond.ipp.bord} provide a set of measures which are admissible. The admissible measures are determined through the requirement that $\sum_j g^{ij}\partial_j\log \rho = A_i$, with $\deg_a(A_i)\leq a_i$, or in other terms $\Gamma(\log \rho, x_i)= A_i$.  When the reduced equation of the boundary is $\{\det(g)=0\}$, then we have described all the measures in Theorem~\ref{main.thm.gal}.  But when some factor of $\det(g)$ does not appear in the reduced equation of $\partial\Omega$, it is nor excluded that those factor may provide some other admissible measure (see~\cite{BOZ}). However, in  dimension two and for the usual degree, where we are able do provide a complete description of all possible models, this situation never appears and we wonder if this  may appear at all. 

The fact that the boundary is included into $\{\det(g)\}=0$ allows to restrict in general to one of the connected components of the complement of this set, so that the metric may never degenerate inside $\Omega$. But it may happen (although we have no example for this) that there exist some solutions of this problem for which the solution $(g^{ij})$ is not positive inside any bounded  region bounded by $\{Q=0\}$.

But the determination of all possible admissible boundaries (that is the curves for which equation~\eqref{cond.ipp.bord} have a non trivial solution) is  a much harder problem. The possibility for an algebraic surface to have a non trivial solution in equation~\eqref{cond.ipp.bord} is a very strong requirement, as we shall see next, and this reduces considerably the  possible models.

\section{Classification with usual degree in dimension 2\label{sec.dim2}}

In this Section, we reduce to the two dimensional case, and moreover to  the usual degree $a_1=a_2=1$. In this situation, the problem is then invariant under affine transformation, and this allows to use classical tools of algebraic geometry to reduce the problem. The coefficients $(g^{ij})$ have at most  degree $2$ and the boundary maximum degree $4$.
The main result is the following

\bthm In dimension $2$ and for the usual degree $a_1=a_2=1$,  and  up to affine transformations,   there exist exactly 11 bounded sets $\Omega$ (with piecewise $\cC^1$ boundary)  corresponding to a polynomial model. Their boundaries are (see pictures in Section~\ref{sec.pictures}):
the triangle, the circle, the square, two coaxial parabolas, a parabola, a tangent line and a line parallel to the axis, a parabola and two tangent lines, a cuspidal cubic and a tangent line, a cuspidal cubic and a line passing through the infinite point   of the cubic, a nodal cubic, the swallow tail and the deltoid curve.

In all the models, the only possible values for the measure are the one described in Theorem~\ref{main.thm.gal}. When the boundary has maximal degree, then the metric $(g^{ij})$ is unique up to scaling, and correspond to a constant curvature metric, either $0$ or $1$ (after appropriate scaling). 

There are models (triangle, circle) where the metric $(g^{ij})$ is not unique.

\ethm

\brmq The previous assertion is not completely exact. The family described by two axial parabolas  are not reducible one to the other under affine transformations. But a simple quadratic transformation does the job.

\ermq

\bpf It is out of scope to give the proof here, which is quite lengthy and technical. But it relies on some simple argument. The main point is to show (using appropriate change of coordinates allowed by the affine invariance of the problem) that there may be no flex point and no flat point on the boundary. That is  (in complex variables), that one may no find an analytic branch locally of the form $y= x^3+ o(x^3)$ or $y= x^4+ o(x^4)$. This is done through the local study  of equation~\eqref{cond.ipp.bord}. Such points correspond to singular points of the dual curve. But there is a balance between  the singular points of the curve, of it's dual curve, and the genus of the curve (seen as a compact Riemann surface), known as Plucker's formulae. This allows to show that $\partial\Omega$ must indeed have many singular points, which list is easy to write since the degree is low (here 4). It remains to analyze all the possible remaining cases. See~\cite{BOZ} for details. 
\epf

Observe that the triangle and the circle case where already described in Section~\ref{sec.spheres} as image of the two dimensional sphere. But even in this case, equation~\eqref{cond.ipp.bord} produces other metrics $(g^{ij})$ than the one already described. If one considers a single entry of a $SU(d)$ matrix, then it corresponds to a polynomial model in the unit disk which is one of these exotic metrics on the unit ball in $\bbR^2$.  Typically, on the circle, one may add to the generator $a(x\partial_y-y\partial_x)^2$, which satisfies the boundary condition and corresponds to some extra  random rotation in the associated stochastic process.  Indeed, all these $11$ models may be described, at least for some half-integer values of the parameters appearing in the description of the measure, as the image of the above mentioned models constructed on spheres, $SO(d)$ or $SU(d)$. But there are also, for other values of the measure, some constructions provided by more sophisticated geometric models, in particular root systems in Euclidean spaces. We refer to \cite{BOZ} for a complete description of theses models.

\section{ Other models with other degrees in dimension 2\label{sec.other.degrees}}

When $\Omega$ is not bounded, equation~\eqref{cond.ipp.bord} is not fully justified. If one restricts our attention to the usual degree and to those boundaries which satisfy this condition, then we obtain only the products of the various one dimensional models, and two extra models which are bounded by a cuspidal cubic or a parabola. In this situation, there are some exponential factors in the measures, as happens in the Laguerre case. When there is no boundary at all, it may be proved (although not easily) that the only admissible measures are the Gaussian ones : they correspond to the product of Ornstein-Uhlenbeck operators, but as is the case with the circle, one may add to the metric some rotational term $(x\partial_y-y\partial_x)^2$, which produces new families of orthogonal polynomials.

Beyond this, one may exhibit some examples  on various bounded sets $\Omega$ with weighted degrees. There is no complete classification in general for such general models at the moment. The reason is that affine invariance is then lost (this is counterbalanced by the fact that  some other polynomial change of variables are allowed), but the local analysis made above is no longer valid. To show how rich this new family of models may be, we just present here some examples.

On $SU(3)$, let $Z$ be the trace of the matrix, considered as a function $Z : SU(3)\mapsto \bbC= \bbR^2$. Then thanks to the fact that there are 3 eigenvalues belonging to the unit circle and product $1$, the characteristic polynomial $\det(M-\Id)$ of an $SU(3)$ matrix may be written as $-X^3+ZX^2- \bar Z X+1$ such that $Z$ itself encodes  the spectral measure.   Applying formulae~\eqref{spect.su3}, and up to a change of $Z$ into $Z/3$ and scaling, one gets
\beq\label{matrix.complex.deltoid}
  \begin{cases}\Gamma(Z,Z)=\bar Z- Z^2,\\
\Gamma(\bar Z, Z)= 1/2(1-Z\bar Z),\\
\Gamma(\bar Z, \bar Z)= Z- \bar Z^2,\\
 \LL Z=  -4Z , \LL \bar Z= -4\bar Z.\end{cases}
 \eeq
 This corresponds indeed with the deltoid model appearing in Section~\ref{sec.dim2}. From these formulae, 
 one sees that functions $F(Z, \bar Z)$ which are symmetric in $(Z,\bar Z)$ are preserved by the image operator. Setting $S= Z+\bar Z$ and $P= Z\bar Z$, this leads to the following polynomial model with degree $\deg_S + 2\deg_P$, with
 $$\bcas
\Gamma(S,S)= 1+S+P-S^2,\\
\Gamma(S,P)= \frac{1}{2} S-2P +S^2-\frac{3}{2} SP,\\
\Gamma(P,P)= P-3SP-3P^2+S^3,\\
\LL S= -4S, \LL (P)= 1-9P.
\ecas$$
Up to some constant, the determinant of the metric is $(4P-S^2)(4S^3-3P^2-12SP-6P+1)$, and the boundary is the domain which is delimited by a parabola and a cuspidal cubic (a degree 5 curve), which are bi-tangent at their intersection point. This leads to a two-parameter family of measures and associated orthogonal polynomials.

One may also construct  more models using discrete symmetry groups. Here are some examples.

We cut the $2$-d sphere into $n$ vertical slices along the meridians, and and look for a basis of functions invariant under the reflections around these meridians. Writing the sphere as  $x_1^2+x_2^2+ x_3^2=1$,  we chose $X= x_3$ and writing in complex notations $x_1+{\bf i}x_2= z$,we  chose $Y= \Re(z^n)$. In polynomials terms $Y= (x_1^2+x_2^2)^{n/2} P_n(\frac{x_1}{\sqrt{x_1^2+x_2^2}})$, where $P_n$ is the $n$-th first kind Chebychef polynomial.  For parity considerations on $P_n$, this is always a polynomial in the variables $(x_1,x_2)$. 
We have  
$$\bcas\Gamma(X,X)= 1-X^2,\\  
\Gamma(X,Y)=  -nXY,\\
\Gamma(Y,Y)= n^2\Big((1-X^2)^{n-1}-Y^2\Big)\\
\LL X-= -2X, ~\LL(Y)= -n(n+1)Y
\ecas
$$
The boundary equation is then $(1-X^2)^n-Y^2=0$, which is irreducible when  $n$ is odd, leading to a one parameter family of measures, and splits into two parts when $n$ is even, leading to a two parameters family.
We may in the previous model look at functions of $(X^2,Y)$ adding a new invariance under symmetry around the hyperplane $\{x_3=0\}$, or of $(X,Y^2)$, or also of $(X^2,Y^2)$ leading to 3 new families.

There are  other ways to construct such  two dimensional models. A general idea is to consider some finite sub-group of $SO(3)$, and extract from the axes of the rotation some subfamily $V_i$ which is invariant under the group action. Then, one considers the homogeneous polynomial $P(x,y,z)
=\prod_i X\cdot V_i$, where $X= (x,y,z)$ and $X\cdot V$ denotes the scalar product. It is also invariant under the group action.  Let $m$ be the degree of this polynomial. With $P$, one constructs a new polynomial $Q= r^{m+1/2}P(\partial_x,\partial_y, \partial_z)r^{-1/2}$, where $r= (x^2+y^2+z^2)^{1/2}$. $Q$ is still homogeneous with degree $m$, invariant under the group action, and moreover harmonic. Then, one looks for the system $X= Q, Y= \Gamma(Q,Q)$ where $\Gamma$ is the spherical Laplace operator on the unit sphere $\bbS^2\subset \bbR^3$. The action of the spherical operator on the pair $(X,Y)$ may lead to a polynomial system. This is not always the case however. For example, with the symmetry group of the icosahedron, there are 3 homogeneous polynomials  $P_6,P_{10}$ and $P_{15}$, with degrees $6,10$ and $15$ which generate all homogeneous polynomials which are invariant under the group action. The technology that we provided works starting from $P_6$ but not from $P_{15}$. The associated formulae are too complicated to be given here. For the reader interested in explicit computations (see~\cite{b.meyer54} for more examples), the explicit value of $P_6$ is as follows (with $c=(1+\sqrt{5})/2$)
$$P_6(x,y,z)= (c^2x^2-y^2)(c^2y^2-z^2)(c^2z^2-x^2).$$

\section{ Higher dimensional models\label{sec.higher.dim}}

The technology which allows to describe all the bounded dimensional models for the usual degree is not available in higher dimension, mainly because of the lack of the analogues of Plucker's formulae. The many models issued from Lie group action that we produced so far provide many families, with various degrees. Beyond these explicit constructionss, and sticking to the bounded models with the usual degree in dimension 2, it is worth to observe that one may produce new admissible sets by double cover. More explicitly, as soon as we have a model in dimension $d$, with reduced boundary equation $P(x)=0$, one may look for models in dimension $d+1$ with boundary equatio $y^2-P(x)=0$ ($y$ being the extra one dimensional variable. It turns out  that this produces a new model for every two dimensional model which have no cusps or tangent lines as singular points (that is, in our setting, the circle, the triangle, the square, the double parabola 
and the nodal cubic).  The boundary has  no longer maximal degree, even if the starting model has, and the metric is not unique in general. Moreover, even in the simplest cases, the curvature is not constant. 

We may then pursue the construction adding new dimensions. The reason why this works (together with the obstruction about singular points) remains mysterious.   Most of the questions regarding these constructions and others remain open at the moment.

\section{Pictures\label{sec.pictures}}
In this Section, we give the various pictures for the 11 bounded models in dimension 2 with natural degree. We give the reduced equation of the boundary, and we indicate when the metric is unique, up to a scaling factor. When it is unique, we indicate the cases when the curvature is constant, and what is it's sign. It is worth to observe that all the models with maximal degree  (here 4) have a unique constant curvature metric.
\begin{figure}[ht!]
\begin{center}
\includegraphics[width=.3\textwidth]{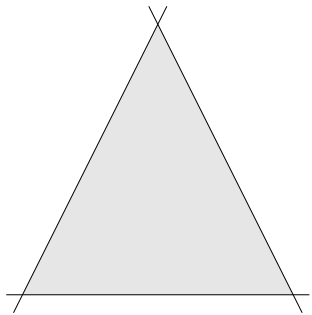}
\caption{Triangle  : $xy(1-x-y)=0$, metric not unique}
\end{center}
\end{figure}
\begin{figure}[ht!]
\begin{center}
\includegraphics[width=.3\textwidth]{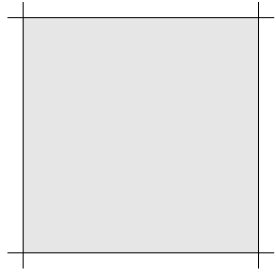}
\caption{Square : $(1-x^2)(1-y^2)=0$, one  metric ,  curvature 0 }
\end{center},
\end{figure}
\begin{figure}[ht!]
\begin{center}
\includegraphics[width=.3\textwidth]{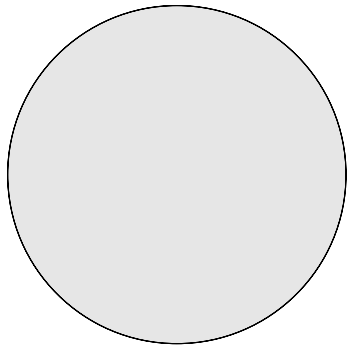}
\caption{Circle  : $x^2+y^2=1$, metric not unique }
\end{center}
\end{figure}
\begin{figure}[ht!]
\begin{center}
\includegraphics[width=.3\textwidth]{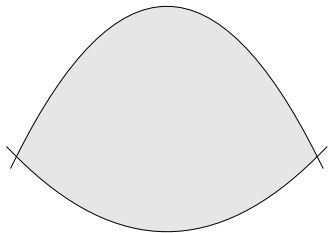}
\caption{Double parabola  : $(y+1-x^2)(y-1+ax^2)=0$, one metric, curvature 1}
\end{center}
\end{figure}
\begin{figure}[ht!]
\begin{center}
\includegraphics[width=.3\textwidth]{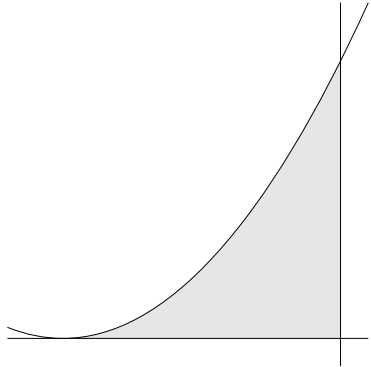}
\caption{Parabola with 1 tangent and one secant line  : $(y-x^2)y(x-1)=0$, one metric, curvature 1}
\end{center}
\end{figure}
\begin{figure}[ht!]
\begin{center}
\includegraphics[width=.3\textwidth]{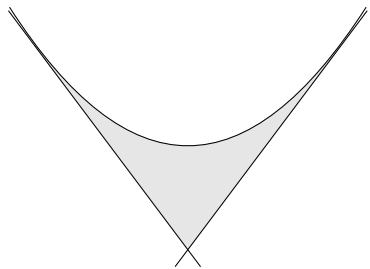}
\caption{Parabola with 2 tangents   : $(y-x^2)(y+1-2x)(y+1+2x)=0$, one metric, curvature 0}
\end{center}
\end{figure}
\begin{figure}[ht!]
\begin{center}
\includegraphics[width=.3\textwidth]{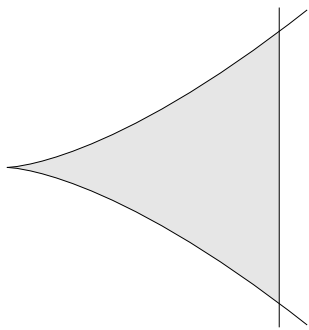}
\caption{Cuspidal cubic with secant line   : $(y^2-x^3)(x-1)=0$, one metric, curvature 1}
\end{center}
\end{figure}
\pagebreak
\begin{figure}[ht!]
\begin{center}
\includegraphics[width=.3\textwidth]{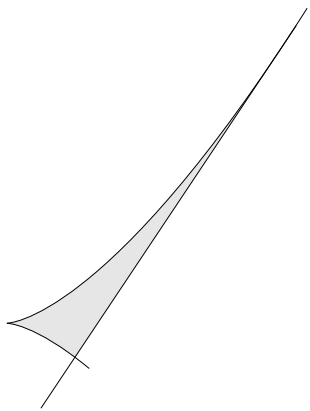}
\caption{Cuspidal cubic with tangent line   : $(y^2-x^3)(3x-2y-1)=0$, one metric, curvature 1}
\end{center}
\end{figure}
\begin{figure}[ht!]
\begin{center}
\includegraphics[width=.3\textwidth]{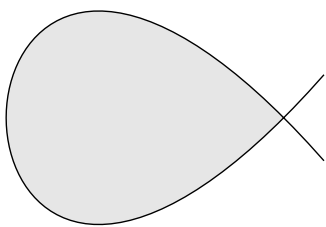}
\caption{Nodal cubic   : $y^2-x^2(1-x)=0$, one metric, non  constant curvature}
\end{center}
\end{figure}
\begin{figure}[ht!]
\begin{center}
\includegraphics[width=.3\textwidth]{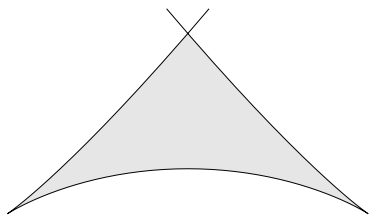}
\caption{Swallow tail   : $4x^{2}-27x^{4}+16y-128y^{2}-144x^{2}y+256y^{3}=0$, one metric, curvature 1}
\end{center}
\end{figure}
\begin{figure}[ht!]
\begin{center}
\includegraphics[width=.3\textwidth]{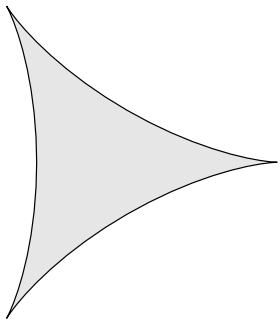}
\caption{Deltoid  : $(x^{2}+y^{2})^{2}+ 18(x^{2}+y^{2})- 8x^{3}+24xy^{2}-27=0$, one metric, curvature 0}
\end{center}
\end{figure}

\bibliographystyle{amsplain}   

\providecommand{\bysame}{\leavevmode\hbox to3em{\hrulefill}\thinspace}
\providecommand{\MR}{\relax\ifhmode\unskip\space\fi MR }
\providecommand{\MRhref}[2]{%
  \href{http://www.ams.org/mathscinet-getitem?mr=#1}{#2}
}
\providecommand{\href}[2]{#2}
\begin{thebibliography}{}

\end{thebibliography}


\begin{thebibliography}{10}

\bibitem{ArCoVa01}
J.~Arves\`u, J.~Coussement, and W.~Van~Assche, \emph{Some discrete multiple
  orthogonal polynomials},
  Preprint,www.wis.kuleuven.ac.be/applied/intas/RomeJon.pdf, 2001.

\bibitem{Askey74}
R.~Askey, \emph{Jacobi polynomials. {I}. {N}ew proofs of {K}oornwinder's
  {L}aplace type integral representation and {B}ateman's bilinear sum}, SIAM J.
  Math. Anal. \textbf{5} (1974), 119--124. \MR{MR0385197 (52 \#6062)}

\bibitem{BGL}
D.~Bakry, I.~Gentil, and M.~Ledoux, \emph{Analysis and {G}eometry of {M}arkov
  {D}iffusion {O}perators}, Grund. Math. Wiss., vol. 348, Springer, Berlin,
  2013.

\bibitem{BakM}
D.~Bakry and O.~Mazet, \emph{Characterization of markov semigroups on {$\mathbb
  R$} associated to some families of orthogonal polynomials.}, S{\'e}minaire de
  Probabilit{\'e}s XXXVII, vol. Lecture Notes in Math. 1832, Springer, Berlin,
  2003., pp.~60--80.

\bibitem{BOZ}
D.~Bakry, S.~Orevkov, and M.~Zani, \emph{Orthogonal polynomials and diffusions
  operators}, preprint, 2013.

\bibitem{BraakMeul}
B.L.J. Braaksma and B.~Meulenbed, \emph{Jacobi polynomials as spherical
  harmonics}, Nederl. Akad. Wetensch. Proc. Ser. A 71=Indag. Math. \textbf{30}
  (1968), 384--389.

\bibitem{BracialiPerezPinar}
Cleonice~F. Bracciali, Teresa~E. P{{\'e}}rez, and Miguel~A. Pi{\~n}ar,
  \emph{Stieltjes functions and discrete classical orthogonal polynomials},
  Comput. Appl. Math. \textbf{32} (2013), no.~3, 537--547. \MR{3120139}

\bibitem{ZerBrink}
H.C. Brinkman and F.~Zernike, \emph{Hypersph{\"a}rische funkktionen und die in
  sph{\"a}rischen bereichen orthogonalen polynome}, Nederl. Akad. Wetensch.
  Proc \textbf{38} (1935), 161--173.

\bibitem{Chihara78}
T.~S. Chihara, \emph{An introduction to orthogonal polynomials}, Gordon and
  Breach Science Publishers, New York, 1978, Mathematics and its Applications,
  Vol. 13. \MR{0481884 (58 \#1979)}

\bibitem{DamMuYes13}
David Damanik, Paul Munger, and William~N. Yessen, \emph{Orthogonal
  {P}olynomials on the {U}nit {C}ircle with {F}ibonacci {V}erblunsky
  {C}oefficients, {II}. {A}pplications}, J. Stat. Phys. \textbf{153} (2013),
  no.~2, 339--362. \MR{3101200}

\bibitem{davies}
E.~B. Davies, \emph{Heat kernels and spectral theory}, Cambridge Tracts in
  Mathematics, vol.~92, Cambridge University Press, Cambridge, 1990.
  \MR{MR1103113 (92a:35035)}

\bibitem{DX}
C.~Dunkl and Y.~Xu, \emph{Orthogonal polynomials of several variables.},
  Encyclopedia of Mathematics and its Applications, vol.~81, Cambridge
  University Press, Cambridge, 2001.

\bibitem{Heck}
G.J. Heckman, \emph{Root sytems and hypergeometric functions i}, Compositio
  Math. \textbf{64} (1987), 353--373.

\bibitem{Heck1}
G.J. Heckmann, \emph{Dunkl operators}, S\'eminaire Bourbaki 828, 1996--97, vol.
  Ast\'erisque, SMF, 1997, pp.~223--246.

\bibitem{KoornWind}
T.~Koornwinder, \emph{Explicit formulas for special functions related to
  symmetric spaces}, Proc. Symp. Pure Math \textbf{26} (1973), 351--354.

\bibitem{Koorn1}
\bysame, \emph{Orthogonal polynomials in two variables which are eigenfunctions
  of two algebraically independent partial differential operators. i.}, Nederl.
  Akad. Wetensch. Proc. Ser. A 77=Indag. Math. \textbf{36} (1974), 48--58.

\bibitem{Koorn2}
\bysame, \emph{Orthogonal polynomials in two variables which are eigenfunctions
  of two algebraically independent partial differential operators. ii.},
  Nederl. Akad. Wetensch. Proc. Ser. A 77=Indag. Math. \textbf{36} (1974),
  59--66.

\bibitem{Koorn3}
\bysame, \emph{Orthogonal polynomials in two variables which are eigenfunctions
  of two algebraically independent partial differential operators. iii.},
  Nederl. Akad. Wetensch. Proc. Ser. A 77=Indag. Math. \textbf{36} (1974),
  357--369.

\bibitem{Koorn4}
\bysame, \emph{Orthogonal polynomials in two variables which are eigenfunctions
  of two algebraically independent partial differential operators. iv.},
  Nederl. Akad. Wetensch. Proc. Ser. A 77=Indag. Math. \textbf{36} (1974),
  370--381.

\bibitem{KrallS}
H.L. Krall and I.M. Sheffer, \emph{Orthogonal polynomials in two variables},
  Ann. Mat. Pura Appl. \textbf{76} (1967), 325--376.

\bibitem{Macdo}
I.~G. Macdonald, \emph{Symmetric functions and orthogonal polynomials.},
  University Lecture Series, vol.~12, American Mathematical Society,
  Providence, RI, 1998.

\bibitem{Macdo1}
\bysame, \emph{Orthogonal polynomials associated with root systems.},
  S{\'e}minaire Lotharingien de Combinatoire, vol.~45, Universit{\'e} Louis
  Pasteur, Strasbourg, 2000.

\bibitem{Maze}
O.~Mazet, \emph{Classification des semi--groupes de diffusion sur {$\mathbb R$}
  associ{\'e}s {\`a} une famille de polyn{\^o}mes orthogonaux}, S{\'e}minaire
  de probabilit{\'e}s XXXI Lecture Notes in Mathematics (J.~Az{\'e}ma et~al,
  ed.), vol. 1655, Springer--Verlag, 1997, pp.~40--53.

\bibitem{b.meyer54}
Burnett Meyer, \emph{On the symmetries of spherical harmonics}, Canadian J.
  Math. \textbf{6} (1954), 135--157. \MR{0059406 (15,525b)}

\bibitem{Sherman}
T.O. Sherman, \emph{The helgason fourier transform for compact riemannian
  symmetric spaces of rank one.}, Acta Mathematica \textbf{164} (1990),
  no.~1-2, 73--144.

\bibitem{szego75}
G.~Szeg{\H{o}}, \emph{Orthogonal polynomials}, fourth ed., American
  Mathematical Society, Providence, R.I., 1975, American Mathematical Society,
  Colloquium Publications, Vol. XXIII.

\bibitem{varopouloscoulhonsaloffcoste92}
N.~Th. Varopoulos, L.~Saloff-Coste, and T.~Coulhon, \emph{Analysis and geometry
  on groups}, Cambridge Tracts in Mathematics, vol. 100, Cambridge University
  Press, Cambridge, 1992.

\end{thebibliography}
\providecommand{\bysame}{\leavevmode\hbox to3em{\hrulefill}\thinspace}
\providecommand{\MR}{\relax\ifhmode\unskip\space\fi MR }
\providecommand{\MRhref}[2]{%
  \href{http://www.ams.org/mathscinet-getitem?mr=#1}{#2}
}
\providecommand{\href}[2]{#2}

\end{document}